\documentclass [12pt] {article}
\usepackage[dvips]{graphicx}
\usepackage{amsmath, amssymb}
\usepackage{amssymb,pb-diagram}
 \usepackage{amscd}

\newcommand{\ZZ}{\mathbb Z}
\newcommand{\PP}{\mathbb P}
\newcommand{\QQ}{\mathbb Q}
\newcommand{\CC}{\mathbb C}

\newcommand{\gtS}{{\mathfrak S}}

\newcommand{\gtA}{{\mathfrak A}}

\newcommand{\Gal}{\mathop {\rm Gal}\nolimits}

\newcommand{\GL}{\mathop {\rm GL}\nolimits}

\newcommand{\PGL}{\mathop {\rm PGL}\nolimits}

\newcommand{\Fix}{\mathop {\rm Fix}\nolimits}

\newcommand{\Pic}{\mathop {\rm Pic}\nolimits}

\newcommand{\Aut}{\mathop {\rm {Aut}}\nolimits}

\newcommand{\Cr}{\mathop{\rm {Cr}}\nolimits}

\newcommand{\dom}{\mathop {{\rm dom}}\nolimits}

\newtheorem{thm}{Theorem}[section]
\newtheorem{cor}{Corollary}[section]

\newtheorem{lem}{Lemma}[section]
\newtheorem{defin}{Definition}[section]
\newtheorem{exmple}{Example}[section]
\newtheorem{rem}{Remark}[section]
\newtheorem{qz}{Question}[section]

\newcommand{\proof}{\noindent{\textsl {Proof}.}\hskip 3pt}

\renewcommand{\thesubparagraph}{\theparagraph.\@arabic\c@subparagraph}
\begin{document}

\title{
 A note on embeddings of $S_4$ and $A_5$ into the Cremona group 
 and versal Galois covers
 }

\author{Shinzo BANNAI and  Hiro-o TOKUNAGA\footnote{Key Words and Phrases: versal Galois cover, Cremona embedding, 2000 MSC 14E20, 14L30.}}
\date{\empty}
\maketitle

\begin{abstract} 
In this article, we prove that two versal Galois covers for $S_4$ and $A_5$ introduced in \cite{tokunaga8},
  \cite{tokunaga9} and \cite{tsuchi} are birationally distinct to each other. As a corollary, we obtain
  two non-conjugate embeddings of $S_4$ and $A_5$ into $\Cr_2(\CC)$. 
\end{abstract}

\large
\textbf{Introduction}\label{intro}
\normalsize
\par\bigskip

 Let $X$ and $Y$ be normal projective varieties defined over $\CC$,
the field of complex numbers. A finite surjective
morphism 
$\pi : X \to Y$ is called Galois, if the induced field extension $\CC(X)/\CC(Y)$ of the field of rational
functions is Galois. Given a finite group $G$, we simply call $\pi : X \to Y$ a $G$-cover if it is Galois and 
$\Gal(\CC(X)/\CC(Y)) \cong G$. In \cite{tokunaga8} and \cite{tsuchi}, a notion called ``versal Galois covers"
is introduced, of which definition is as follows:

\begin{defin} \label{def:versal} { \rm  Let $G$ be a finite group. A $G$-cover $\varpi : X \to Y$ is
called a versal Galois cover for $G$ or a versal $G$-cover if it satisfies
the following property:

\medskip

For any $G$-cover $\pi: W \to Z$, there exists a $G$-equivariant rational map $\mu : W \dasharrow X$
 such that

\[
\mu(W) \not\subset \Fix(X, G),
\]
where $\Fix(X, G):=\{x \in X \mid \mbox{the stabilizer group at $x$, $G_x \neq \{1\}$}\}$.
}
\end{defin}

\begin{rem} {\rm The rational map $\mu$ induces a rational map $\overline{\mu} : Z \dasharrow Y$. 
Concerning this
rational map $\overline{\mu}$, 
there exists a Zariski open set $U$ such that $(i)$ $U \subset \dom(\overline\mu)$, $\dom(\bullet)$ being the domain 
of a rational
map $\bullet$, and $(ii)$ $\pi^{-1}(U)$ is birationally equivalent to $U\times_YX$ over $U$. (see \cite{tokunaga9},
Proposition 1.2).
}
\end{rem}

The notion of versal $G$-covers implicitly appeared in \cite{namba1} and \cite{namba2}  as the ``pull-back" 
construction of $G$-covers, and Namba showed that there exists a versal $G$-cover of dimension $\sharp(G)$ for
any finite group $G$. Namba's model, however, has too large dimension for practical use.

By Theorem 6.2 in \cite{buhler-reichstein}, there exists a $1$-dimensional versal $G$-cover if and only if
$G$ is either a cyclic group or a dihedral group of order $2n$ ($n$: odd).
As an next step, 
in \cite{tokunaga8}, \cite{tokunaga9} and \cite{tsuchi}, we studied $2$-dimensional versal $G$-covers and gave
some explicit examples. 
 For a finite subgroup $G$ in  $\GL(n, \ZZ)$, Bannai and  Tsuchihashi construct versal $G$-covers of dimension $n$ by using toric geometry in \cite{bannai} and \cite{tsuchi}.
Among explicit examples in  \cite{tokunaga8}, \cite{tokunaga9}, 
two different versal $G$-covers, $\varpi_{G, 1} : X_1 \to Y_1$ and
$\varpi_{G, 2} : X_2 \to Y_2$ are given for the case when $G$ is  $S_4$, the symmetric group of $4$-letters
and $A_5$, the alternating group of $5$-letters (see \S 1 for description of $X_1$ and $X_2$). By definition of 
versal $G$-covers, there exist $G$-equivariant rational maps $\mu_1 : X_1 \dasharrow X_2$ and
$\mu_2 : X_2 \dasharrow X_1$ such that $\mu_1(X_1) \not\subset \Fix (X_2, G)$ and $\mu_2(X_2) \not\subset
\Fix(X_1, G)$. Moreover, $X_1$ and $X_2$ are del-Pezzo surfaces. Under these circumstances, it may be natural
to raise a question as follows:

\begin{qz} \label{qz:main} {\rm Let $G$ be either $S_4$ or $A_5$. Let $\varpi_{G, 1} : X_1 \to Y_1$ and
$\varpi_{G, 2} : X_2 \to Y_2$ be versal $G$-covers as above. Does there exist any $G$-equivariant
birational map from $X_1$ to $X_2$?
}
\end{qz}

In this note, we consider Question \ref{qz:main} and prove the following:

\begin{thm}\label{thm:main} {There exists no $G$-equivariant birational map from $X_1$ to $X_2$}
\end{thm}

Since both $X_1$ and $X_2$ are rational, their birational automorphism group is the $2$-dimensional Cremona
group $\Cr_2(\CC)$.  For $G = S_4, A_5$, we have two different embeddings $\eta_i : G \to \Cr_2(\CC)$ $(i = 1, 2)$
via $G \subset \Aut(X_i) \subset \Cr_2(\CC) (i = 1, 2)$. Our theorem implies that $\eta_1(G)$ is \textit{not}
conjugate to $\eta_2(G)$ in $\Cr_2(\CC)$. Combining Proposition 0.3 $(i)$ in \cite{tokunaga9}, we have the following
corollary:

\begin{cor}\label{cor:main}{Both $S_4$ and $A_5$ have at least $3$ non-conjugate embedding into $\Cr_2(\CC)$.}
\end{cor}

Our results could be found in old literatures such as \cite{kantor} and \cite{wiman}, but we would like to
emphasize that our question comes from the study of versal $G$-covers, which is a rather new notion. Also conjugate classes of finite subgroups of $\Cr_2(\CC)$ 
have been studied by several mathematicians (\cite{bayle-beauville}, \cite{beauville1}, \cite{beauville2}, \cite{defernex}, \cite{iskovskikh1}).  The notion of versal $G$-covers may add another interest
to this subject.

This article goes as follows. We first give a detailed description of the versal $G$-covers $\varpi_{G,i} :
X_i \to Y_i$ $(i = 1, 2)$ in \S 1. In \S 2, we explain our main tool, ``Noether's inequality," which plays important role in \cite{iskovskikh1} and \cite{iskovskikh2}. 
We prove Theorem \ref{thm:main} in \S 3. In \S 4, we consider rational maps between $X_1$ and $X_2$ in the
case of $G = S_4$.

\bigskip

\textsl{Acknowledgement} A key step of this note was done during the second author's visit to Ruhr Universit\"at
Bochum. He thanks Professor A. Huckleberry for his comments and hospitality.


\section{Versal $S_4$- and $A_5$-covers: Two examples}

\subsection{Versal $S_4$-covers}\label{S_4covers}

Let $S_4$ be the symmetric group of $4$-letters. In order to describe $S_4$, we use
the following representations:

\par\medskip

\begin{tabular}{ll} 
Generators: & $\sigma, \tau, \lambda_1, \lambda_2$ \\

Relations:  & $\sigma^2 = \tau^3 = \lambda_1^2 = \lambda_2^2 = 1$ \\
  & $\sigma\tau = \tau^2\sigma$, $\lambda_1\lambda_2 = \lambda_2\lambda_1$ \\
  & $\sigma\lambda_1\sigma = \lambda_1\lambda_2$, $\sigma\lambda_2\sigma =\lambda_2$ \\
  & $\tau^2\lambda_1\tau = \lambda_1\lambda_2$, $\tau^2\lambda_2\tau = \lambda_1$
\end{tabular}

Let $\rho: S_4 \to \GL(3,\CC)$ be a faithful irreducible representation as follows:
\[
\begin{array}{cc}
 \sigma  \mapsto \left( \begin{array}{ccc}
                        0 & 1 & 0 \\
                        1 & 0 & 0 \\
                        0 & 0 & 1
                         \end{array} \right ), &  \tau \mapsto  \left ( \begin{array}{ccc}
                                                                     0 & 0 & 1 \\
                                                                     1 & 0 & 0 \\
                                                                     0 & 1 & 0 
                                                                      \end{array} \right ), \\
                            & \\                                                             
\lambda_1   \mapsto \left( \begin{array}{ccc}
                        -1 & 0 & 0 \\
                        0 & 1 & 0 \\
                        0 & 0 & -1
                         \end{array} \right ), & \lambda_2 \mapsto  \left ( \begin{array}{ccc}
                                                                     -1& 0 &0 \\
                                                                     0& -1 &0 \\
                                                                     0& 0 & 1 
                                                                      \end{array} \right ).     
                                                                      \end{array}                                                          
\]

\par\medskip

\textbf{Versal $S_4$-cover $\varpi_{S_4, 1} : X_1 \to Y_1$}

\par\medskip

Let $X_1$ be a surface in $\PP^1\times\PP^1\times\PP^1$ defined by the equation
\[
x_0y_0z_0 - x_1y_1z_1 = 0,
\]
where $([x_0,x_1], [y_0,y_1], [z_0,z_1])$ denotes the homogeneous coordinates. Put $x =x_1/x_0, y= y_1/y_0,
z = z_1/z_0$. Define an $S_4$-action on $\PP^1\times\PP^1\times\PP^1$ as follows:
\begin{eqnarray*}
(x, y, z)^{\sigma} &=& (y, x, z),\\
(x, y, z)^{\tau} &=& (y, z, x),\\
(x, y, z)^{\lambda_1} &=& (-x, y, -z),\\
(x, y, z)^{\lambda_2} &=& (-x, -y, z).
\end{eqnarray*}
The defining equation of $X_1$ is invariant under this $S_4$-action. Hence $S_4$ acts on $X_1$. Put $Y_1 = X_1/G$ and denote the quotient morphism by $\varpi_{S_4, 1} : X_1 \to Y_1$. By \cite{tokunaga8} and \cite{tsuchi},  
$\varpi_{S_4, 1} : X_1 \to Y_1$ is a versal $S_4$-cover.

We look into some properties of $X_1$ with respect to this $S_4$-action for later use. We first remark that $X_1$ is a del-Pezzo surface of degree $6$, i.e., $X_1$ is obtained by blowing-up at distinct $3$ points of $\PP^2$.

\begin{lem}\label{lem:1-1}{ The divisor given by $x_0y_0z_0 = 0$ is a cycle of rational curves $C_1, C_2,\dots, C_6$. Each $C_i$ is a smooth rational curve with $C_i^2 = -1$.}
\end{lem}

\proof Let $p_{12} : \PP^1 \times \PP^1\  \times \PP^1 \to \PP^1\times \PP^1$ be the projection to the first two factors.
By its defining equation, we infer that the restriction of $p_{12}$ to $X_1$ is  the blowing-up of $\PP^1\times\PP^1$ at
$([1, 0], [0,1])$ and  $([0,1], [1,0])$. Our statement easily follows from this observation.

\par\bigskip

\begin{lem}\label{lem:1-2}{ Let $\Pic(X_1)$ be the Picard group of $X_1$. Then the $S_{4}$ invariant part $\Pic^{S_4}(X_1) = \ZZ(-K_{X_1})$}.
\end{lem}

\proof $-K_{X_1} \sim \sum_{i=1}^6 C_i$, and one can easily check that the divisor class in the right hand generates 
$\Pic^{S_4}(X_1)$.

\par\bigskip

For $x \in X_1$, we put $d_x = \sharp O(x)$, where $O(x)$ denotes the orbit of $x$. For later use, we study points with $d_x< 6$.

\begin{lem}\label{lem:1-3}{ $(i)$ There are no  points with $d_x = 1, 2, 5$.

$(ii)$ There are $12$ points with $d_x = 4$ as follows:
\[
\begin{array} {cccc}
 R_{11} (1,1,1), & R_{12} (1, -1, -1), & R_{13} (-1, -1, 1), & R_{14} (-1, 1, -1), \\
 R_{21} (\omega, \omega, \omega), & R_{22} (\omega, -\omega, -\omega), & R_{23} (-\omega, -\omega, \omega), & R_{24} (-\omega, \omega, -\omega),   \\
 R_{31} (\omega^2, \omega^2, \omega^2), & R_{32} (\omega^2, -\omega^2, -\omega^2), & R_{33}  (-\omega^2, -\omega^2, \omega^2), &
 R_{34} (-\omega^2, \omega^2, -\omega^2), 
 \end{array}
 \]
 where the coordinates means  the affine coordinates $(x, y, z)$ and $\omega = \exp (2\pi\sqrt{-1}/3)$.
 These 12 points consist of $3$-orbits.
 
 $(iii)$ There are $6$ points with $d_x = 3$ as follows:
 \[
\begin{array} {ccc}
 P_1 ([0,1], [1, 0], [0,1]), & P_2 ([1,0], [0,1], [0,1]), & P_{3} ([0,1], [0,1], [1, 0]), \\
 Q_1([1,0,], [1,0], [0,1]), & Q_2([1,0], [0,1], [1,0]), & Q_3 ([0,1], [1,0], [1,0]).  
 \end{array}
 \]
 These $6$ points consist of $2$ orbits.
 }
 \end{lem}
 
 \proof Note that $\tau$ acts on the divisor $x_0y_0z_0 = 0$ freely and the subgroup $\langle \lambda_1, \lambda_2 \rangle$ has no fixed points
 on the affine surface $xyz = 1$.  Taking these observation into account, we can easily check the above statement by direct computation.
 
 \begin{lem}\label{lem:1-4}{Divisors given by the equations $x_1 = \omega^i x_0$  $(i = 0, 1, 2)$ are rational curves with self-intersection number $0$.}
 \end{lem}
 \proof By the proof of Lemma \ref{lem:1-1}, we infer that the divisors as above come from those in $\PP^1\times \PP^1$ with self-intersection number $0$ and
 all of these divisors in $\PP^1\times \PP^1$ do not pass through $([1, 0], [0,1])$ and  $([0,1], [1,0])$. This implies our statement.
 
 \bigskip

 \textbf{Versal $S_4$-cover $\varpi_{S_4, 2} : X_2 \to Y_2$}
 
 \medskip
 
 By Proposition 4.1 (ii) in \cite{tokunaga8}, we have a versal $S_4$-cover $\PP^2 \to \PP^2/S_4$ by using $\rho$. Put $X_2 = \PP_2$, $Y_2 = \PP^2/S_4$ and
 $\varpi_{S_4, 2}=$ the quotient morphism.

 \subsection{Versal $A_5$-covers}
 
 We first start with the following lemma.
 
 \begin{lem}\label{lem:1-5}{ Let $S$ be a smooth projective surface such that $A_5$ acts faithfully on $S$.  Let $d_x$ be the one as in Lemma \ref{lem:1-3}. Then there exists  no point $x$ on $S$ with $d_x < 5$.}
 \end{lem}
 
 \proof Case $d_x = 1$. Assume that there exists a point $x$ with $d_x = 1$. Then we have a non-trivial homomorphism $\eta : A_5 \to \GL(T_x S)$, where $T_xS$ is the  tangent plane at $x$. Since $A_5$ is simple, $\eta$ is injective, but $A_5$ has no $2$-dimensional faithful representations.
 
 Case $d_x = 2$. Assume that such a point $x$ exists. The stabilizer group $G_x$ at $x$ is a subgroup of $A_5$ with index $2$. This implies that $A_5 \rhd G_x$,
 which is impossible.
 
 Case $d_x = 3, 4$. Assume that such a point exists. Then we have a non-trivial homomorphism from $A_5$ to the symmetric group of either $3$ or $4$ letters.
 The kernel of this homomorphism is a non-trivial normal subgroup, which is a contradiction.
 
 \bigskip

 \textbf {Versal $A_5$-cover $\varpi_{A_5, 1}: X_1 \to Y_1$}

\medskip
Let $X_1$ be a del-Pezzo surface of degree $5$. It is known that $\Aut(X_1) \cong S_5$ (e.g., see \cite{koit}). Based on the result in \cite{hashimoto-tsunogai},
we have proved $X_1 \to X_1/A_5$ is a versal $A_5$-cover. Put $Y_1 = X_1/A_5$ and $\varpi_{A_5, 1}: X_1 \to Y_1$. Note that $\text{Pic}^{A_{5}}(X_1)=\ZZ(-K_{X})$ in this case also.

\bigskip

\textbf{Versal $A_5$-cover $\varpi_{A_5, 2}: X_2 \to Y_2$}

\medskip

Let $\rho : A_5 \to \GL(3, \CC)$ be a faithful irreducible representation. By using $\rho$, we obtain another versal $A_5$-cover $\varpi_{A_5, 2} : X_2 \to Y_2$
in the same manner as in the versal $S_4$-cover $\varpi_{S_4, 2} : X_2 \to Y_2$, i.e. $X_{2}=\PP^{2}$, $Y_{2}=\PP^{2}/A_{5}$, and $\varpi_{A_5, 2}=$ the quotient morphism.

\section{Noether's inequality}

 In this section we explain Noether's inequality in our setting. The proof is  identical to the proof of the general form of Noether's inequality given in \cite{iskovskikh2}. We only need to keep in mind that we are using $G$-invariant linear systems.
 
 Let $X$ and $X^{\prime}$ be smooth projective surfaces with $G$-action.  Let $\mathcal{K}_{X}$ (resp. $\mathcal{K}_{X^{\prime}}$) be the canonical linear system of $X$ (resp. $X^{\prime}$).  Let $\Phi\colon X \dasharrow X^{\prime}$ be a $G$-equivariant birational map. 
Let $\mathcal{H}_{X^{\prime}}$ be a $G$-invariant variable linear system of divisors on $X^{\prime}$ which does not have any fixed components. Let $\mathcal{H} _{X}=\Phi^{-1}(\mathcal{H}_{X\prime})$ be the proper inverse image of $\mathcal{H}_{X^{\prime}}$. Note that $\Phi$ is $G$-equivariant, so $\mathcal{H}_{X}$ is also $G$-invariant.  
 
Let $\eta\colon X_{N}\rightarrow X$ be the $G$-equivariant resolution of indeterminacies of $\Phi$ \cite{reichstein-youssin}. It is a composition of $G$-equivariant blow-ups along smooth centers, which are blow-ups along 0-dimensional $G$-orbits $O(x)$ in our case. Let $\psi=\Phi\circ\eta$.

\[
\begin{diagram}
\node{\eta\colon X_{N}}\arrow{e,t}{\eta_{N,N-1}}\node{X_{N-1}}\arrow{e,t}{\eta_{N-1,N-2}}\node{\ldots}\arrow{e,t}{\eta_{2,1}}\node[1]{X_{1}}\arrow{e,t}{\eta_{1,0}}\node{X_{0}=X}
\end{diagram}
\]
\[
\begin{diagram}
\node{X_{N}}\arrow{s,l}{\eta}\arrow{se,t}{\psi}\\\node{X}\arrow{e,r,..}{\Phi} 
\node[1]{Y}
\end{diagram}
\]
$\eta_{i+1,i}$ is a blow-up along a 0-dimensional $G$-orbit $O(x_{i})$.
Let $\eta_{j,i}=\eta_{j,j-1}\circ\cdots\circ\eta_{i+1,i}$ ($N\geq j>i+1\geq 1$), $\eta_{N,N}=\text{id}_{X_{N}}$. Let $\mathcal{H}_{X_{N}}$ be the proper transform of $\mathcal{H}^{\prime}$ on $X_{N}$. Let $H_{\bullet}$ and $K_{\bullet}$ be a member of $\mathcal{H}_{\bullet}$ and $\mathcal{K}_{\bullet}$ respectively, where $\bullet=X$, $X_{N}$, and $X^{\prime}$. 
Then we have
$$H_{X_{N}}=\eta^{\ast}H_{X}-\sum_{i=0}^{N-1}r(x_{i})\eta_{N,i+1}^{\ast}(E_{i+1})$$
$$K_{X_{N}}=\eta^{\ast}K_{X}+\sum_{i=0}^{N-1}\eta^{\ast}_{N,i+1}(E_{i+1})$$
where $r(x_{i})$ is the multiplicity of a base point $x_{i}\in O(x_{i})$ (a point in the center $O(x_{i})$ of the blow-up $\eta_{i+1,i}$) of $\mathcal{H}_{X}$, and $E_{i}$ is the exceptional divisor of $\eta_{i,i-1}$. We note that $E_{i}$ is a disjoint union of $(-1)$-curves corresponding to the points in $O(x_{i})$, and $r(x_{i})=r(x_{j})$ if $O(x_{i})=O(x_{j})$ since $\mathcal{H}_{X}$ is $G$-invariant.

 \begin{lem}\label{lem:noether}[Noether's Inequality]
 
(i) Under the notation above, suppose that  $\mathcal{H}_{X^{\prime}}+m\mathcal{K}_{X^{\prime}}=\emptyset$ then there exists a 0-dimensional $G$-orbit $O(x)$ consisting of maximal singularities (i.e. base points of $\mathcal{H}$ with multiplicity $r>m$), or else the adjoint linear system $\mathcal{H}_{X}+m\mathcal{K}_{X}$ is empty on $X$. 
 
(ii) Under the same conditions, if there exists a variable family of curves $\mathcal{C}^{\prime}$ such that   $(H_{X^{\prime}}+mK_{X^{\prime}})C^{\prime}<0$ then there exists a 0-dimensional $G$-orbit of maximal singularities, or else there is a curve $C\subset X$ such that   $(H_{X}+mK_{X})C<0$.
 \end{lem}
 
\begin{proof}
$(i)$We have
\begin{align}
H_{X_{N}}+mK_{X_{N}}&=\eta^{\ast}(H_{X}+mK_{X})+\sum_{i=0}^{N-1}(m-r(x_{i}))\eta^{\ast}_{N,i+1}(E_{i+1}) \label{divisor}\\
\intertext{Then by applying ${\psi_{\ast}}$ to both sides, we have}
H_{X^{\prime}}+mK_{X^{\prime}}&=\psi_{\ast}(H_{X_{N}}+mK_{X_{N}}) \notag \\
&=\psi_{\ast}\eta^{\ast}(H_{X}+mK_{X})+\psi_{\ast}(\sum_{i=0}^{N-1}(m-r(x_{i}))\eta^{\ast}_{N,i+1}(E_{i})) \notag 
\end{align}
Since $\mathcal{H}_{X^{\prime}}+m\mathcal{K}_{X^{\prime}}=\emptyset$ by hypothesis the right hand side cannot be an effective divisor, hence $r(x_{i})>m$ for at least one $i$, or else $\mathcal{H}_{X}+m\mathcal{K}_{X}=\emptyset$. 

$(ii)$$\psi^{\ast}(H_{X^{\prime}}+mK_{X^{\prime}})=(H_{X_{N}}+mK_{X_{N}})+F$ where $F$ is the exceptional divisor of $\psi$. Then $\psi^{\ast}\mathcal{C}^{\prime}F=0$. Then we have
$(H_{X_{N}}+mK_{X_{N}})\psi^{\ast}\mathcal{C}^{\prime}<0$. Suppose that $r(x_{i})\leq m$ for all $i$. Then by intersecting both sides of (\ref{divisor}) by $C\in \psi^{\ast}\mathcal{C^{\prime}}$ we find that ${\eta^{\ast}}(H_{X}+mK_{X})\psi^{\ast}\mathcal{C}^{\prime}<0$. Hence $(H_{X}+mK_{X})\eta_{\ast}\psi^{\ast}\mathcal{C}^{\prime}<0$. A general member $C^{\prime}$ of $\eta_{\ast}\psi^{\ast}\mathcal{C}^{\prime}$ may be reducible but we have $(H_{X}+mK_{X})C<0$ for at least one irreducible component of $C^{\prime}$. 
\end{proof}
 

 \section{Proof of Theorem \ref{thm:main}}

 \subsection{ The case of $S_4$}
 
 Suppose that  there exists an $S_4$-equivariant rational map $\Phi : X_1 \dasharrow X_2 (= \PP^2)$.  Let $\Lambda$ be the complete linear system
 given by the class of line $L$, and let $\Phi^{-1}(\Lambda)$ be the proper image of $\Lambda$.  Since the map $\Phi$ is given by $\Phi^{-1}(\Lambda)$,
 $\Phi^{-1}(\Lambda)$ has no fixed components. Also $\Phi^{-1}(\Lambda)$ is $S_4$-invariant.  Hence  $H \in \Phi^{-1}(\Lambda)$ is linearly equivalent
 to $-aK_{X_1}$ for some $a \ge 1$.  Now apply Lemma \ref{lem:noether} to $\Lambda+ a\mathcal{K}_{X_2}$ and  $\Phi^{-1}(\Lambda) + a\mathcal{K}_{X_1}$.  There 
 exists an $S_4$-orbit $O(x)$ of $x \in X_1$ such that $O(x)$ is contained in the base points $\Phi^{-1}(\Lambda)$ and the multiplicity $r$ at each point in
 $O(x)$ is greater than $a$. As any element in $\Phi^{-1}(\Lambda)$ passes through $O(x)$ with multiplicity $r$, we have $a^2K_{X_1}^2 \ge r^2d$, $d$ being
 $\sharp(O(x))$; and we have $d < K_{X_1}^2 = 6$. Hence $O(x)$ is one of the orbits described in Lemma \ref{lem:1-3}. Following 
 Iskovskhih \cite{iskovskikh2}, we call this orbit \textit{a maximal singularity}.
 
 \begin{lem}\label{lem:3-1}{The orbit $O(x)$ with $d = 4$ can not be a maximal singularity.}
 \end{lem}
 
 \proof Let $E_i$ be the divisor on $X_1$ given by $x_1 = \omega^ix_0$ $(i = 0, 1, 2)$ as in Lemma \ref{lem:1-4}.  Suppose that $O((\omega^i, \omega^i, \omega^i))$
 is a maximal singularity, and let $q : {\hat X}_1 \to X_1$ be the blowing-up at $O((\omega^i, \omega^i, \omega^i))$.  Then the linear system
 $q^\ast(\Phi^{-1}(\Lambda)) - r(R_{i1} + R_{i2} + R_{i3}  + R_{i4})$ does not have any fixed components (we identify $R_{ij}$ $(j = 1, 2,3,4)$ with the exceptional curves).
 Let ${\bar E}_i$ be the proper transform of $E_i$. Then
 \[
 (-aq^*K_{X_1} - r \sum_{j=1}^4R_{ij}){\bar E_i} = 2a - 2r < 0.
 \]
 This means that  ${\bar E}_i$ is a fixed component of $q^{\ast}(\Phi^{-1}(\Lambda)) - r(R_{i1} + R_{i2} +R_{i3} + R_{i4})$.

 \begin{lem}\label{lem:3-2}{The orbit $O(x)$ with $d = 3$ can not be a maximal singularity.}
 \end{lem}
 
 \proof  Suppose that $O(P_1) = \{P_1, P_2, P_3\}$ is  a maximal singularity. We may assume that  the irreducible component $C_1$ in the divisor $x_0y_0z_0 = 0$ passes through $P_1$.   Let $q : {\hat X}_1 \to X_1$ be the blowing-up at $O(P_1)$.  Then the linear system
 $q^\ast(\Phi^{-1}(\Lambda)) - r(P_1 + P_2 + P_3)$ does not have any fixed components (we identify $P_i$ $(j = 1, 2, 3)$ with the exceptional curves).
 Let ${\bar C}_1$ be the proper transform of $C_1$. Then
 \[
 (-aq^*K_{X_1} - r \sum_{j=1}^3P_j){\bar C}_1 = a - r < 0.
 \]
 This means that  ${\bar C}_1$ is a fixed component of   $q^{\ast}(\Phi^{-1}(\Lambda)) - r(P_1 + P_2 + P_3)$. 
 
 \bigskip
 
 By Lemmas \ref{lem:3-1} and \ref{lem:3-2}, Theorem \ref{thm:main} for $S_4$ follows.
 
 \subsection {The case of  $A_5$}  By the same argument as in the previous case, the existence of $\Phi$ implies the existence  of an $A_5$-orbit $O(x)$, $x \in X_1$ with 
 $\sharp (O(x)) < 5$. This contradicts Lemma \ref{lem:1-5}.

 
 \section{A remark for versal $S_4$-covers $\varpi_{S_4, 1} : X_1 \to Y_1$  and $\varpi_{S_4, 2} : X_2 \to Y_2$}
 
 By the definition of versality, there exist $S_4$-equivariant rational maps $\mu_1 : X_1 \dasharrow X_2$ and
 $\mu_2 : X_2 \dasharrow X_1$ such that $\mu_{1}(X_{1)}\not\subset\text{Fix}(X_{2},G)$ and  $\mu_{2}(X_{2})\not\subset\text{Fix}(X_{1},G)$. Note that both of  $\mu_i$ $(i = 1, 2)$ are dominant as there exists no $1$-dimensional
 versal $S_4$-cover.  In this section, we give examples of such $\mu_i$ $(i = 1,2)$ such that
 
\par\medskip

$(i)$ both field extensions $\CC(X_1)/\CC(X_2)$ and $\CC(X_2)/\CC(X_1)$ induced by $\mu_1$ and $\mu_2$, respectively, are
cyclic extension of degree $3$, and

\par\medskip

$(ii)$ the field extension $\CC(X_2)/(\mu_2\circ\mu_1)^*(\CC(X_2))$ is Galois and its Galois group is ismorphic to 
$(\ZZ/3\ZZ)^{\oplus 2}$.

\par\bigskip

Let $G$ be a subgroup of $\PGL (3, \CC)$ generated by the following matrices:
\[
\begin{array}{ccc}
A = \left (\begin{array}{ccc}
                   1 & 0 & 0 \\
                   0 & 0 & 1 \\
                   0 & 1 & 0
                    \end{array}\right ), &
                                                       B = \left (\begin{array}{ccc}
                                                                       0 & 0 & 1 \\
                                                                       1 & 0 & 0 \\
                                                                       0 & 1 & 0
                                                                        \end{array} \right ), & 
                                                                                                                C_1 = \left ( \begin{array}{ccc}
                                                                                                                                    -1 & 0 & 0 \\
                                                                                                                                    0 & -1 & 0 \\
                                                                                                                                    0 & 0 & 1
                                                                                                                                    \end{array} \right ), \\
                                                                                                                                   & & \\
                                                                                                                                    
C_2 = \left (\begin{array}{ccc}
                    1 & 0 & 0 \\
                    0 & -1 & 0 \\
                    0 & 0 & -1  
                     \end{array}\right ), &
                                                         D_1 = \left (\begin{array}{ccc}
                                                                             1 & 0 & 0 \\
                                                                             0 & \omega & 0 \\
                                                                             0 & 0 & 1 \end{array} \right ),  &
                                                                                                               D_2 = \left ( \begin{array}{ccc}
                                                                                                                                   1 & 0 & 0 \\
                                                                                                                                   0 & 1 & 0 \\
                                                                                                                                   0 & 0 & \omega
                                                                                                                                    \end{array} \right ),
\end{array}
\]
where $\omega = \exp(2\pi\sqrt{-1}/3)$.

One can easily check the facts below
\begin{enumerate}

\item The order of $\sharp (G)$ of $G$ is $216$.

\item The subgroup, $H_1$,  generated by $A, B, C_1$ and $C_2$ is isomorphic to $S_4$.

\item The subgroup, $H_2$, generated by $D_1$ and $D_2$ is a normal subgroup of $G$. $G$ is a semi-direct product of $H_1$ and
$H_2$.

\item The subgroup, $H_3$, generated by $D_1D_2^2$ is a normal subgroup of $G$.

\item The representation  of $S_4$ by $A$, $B$, $C_1$, and $C_2$  is conjugate to the representation described in Section \ref{S_4covers}, and is given by $\rho(\tau)(\rho(g))\rho(\tau)^{-1}, g \in S_4$.

\end{enumerate}

Let $[X_0, X_1, X_2]$ be the homogeneous coordinate functions on $\PP^2$. We define the action of $G$ by
\[
[X_0, X_1, X_2]^g := [X_0, X_1, X_2]A_g,
\]
where $A_g$ denotes the $3\times 3$ matrix corresponding $g$. Put $x = X_1/X_0$ and  $y = X_2/X_0$. Then
$\CC(\PP^2) \cong \CC(x, y)$ and the $G$-action with respect to $x$ and $y$ is described as follows:
\[
\begin{array}{ccc} 
 (x, y)^A = (y, x), & (x, y)^B = (y/x, 1/x), & (x, y)^{C_1} = (x, -y), \\
(x, y)^{C_2} = (-x, -y), & (x, y)^{D_1} = (\omega x, y) & (x, y)^{D_2} = (x, \omega y).
\end{array}
\]
  
\begin{lem}\label{lem:generator}{ There exists $u, v \in \CC(\PP^2)^{H_3}$ such that $(i)$ $\CC(\PP^2)^{H_3} =
\CC(u, v)$, the invariant field by $H_3$,  and $(ii)$  the $S_4$-action induced by $H_1$ is described as follows:
\[
\begin{array}{cc}
(u, v)^A = (v, u), & (u, v)^B = (v, 1/(uv)), \\
(u, v)^{C_1} = (-u, v), & (u, v)^{C_2} = (-u, -v).
\end{array}
\]
}
\end{lem}

\begin{proof} Put $u = x^2/y$ and $v = y^2/x$. Then one can easily check $u^{D_1D_2^2} = u$ and $v^{D_1D_2^2} = v$.
Hence $u, v \in \CC(\PP^2)^{H_3}$. We show that $\CC(u, v) = \CC(\PP^2)^{H_3}$.  Let $\theta = y/x$. As
$\theta^{D_1D_2^2} = \omega\theta$, $\theta \not\in \CC(\PP^2)^{H_3}$.  Since $x = u\theta$ and $y =
\theta/v$, we have $\CC(\PP^2) = \CC(u, v)(\theta)$. Moreover, $\theta^3 = v/u \in \CC(u, v)$. Hence
$[\CC(\PP^2):\CC(u, v)] = 3$ and this implies that $\CC(u, v) = \CC(\PP^2)^{H_3}$. The second assertion is 
straightforward by our choice of $u$ and $v$. \end{proof}

We now see that the $S_4$-action on $\CC(u, v)$ is identified with the one on $\CC(X_1)$ as follows:

Let $(x, y, z)$ be the affine coordinate of $\PP^1\times \PP^1 \times \PP^1$ introduced in \S 1. The surface $X_1$ is
given by the equation $xyz = 1$ with respect to this coordinate. Put $x = u, y = v, z = 1/xy$, we can check that the
two $S_4$-actions coincide with each other.

Now consider a sequence of fields: 
\[
\CC(\PP^2) \subset \CC(\PP^2)^{H_3} = \CC(X_1) \subset \CC(\PP^2)^{H_2}.
\]
Corresponding to this sequence, we have  rational maps
\[
\mu_2: \PP^2 \dasharrow X_1 
\]
and
\[
\mu_1: X_1 \dasharrow \PP^2/H_2.
\]

Since the subgroup $H_1$ induces an $S_4$-action on each surface,  both $\mu_1$ and $\mu_2$
are $S_4$-equivariant.  Also the quotient surface $\PP^2/H_2$ is  isomorphic to $\PP^2$ and  the $S_4$-actions on
$\PP^2$ and $\PP^2/H_2$ induced by $H_1$ is identified with the one on $X_2 (= \PP^2)$.  
Finally the field extension $(\mu_2\circ\mu_1)^*(\CC(X_2)) \subset \CC(X_2)$ is a $(\ZZ/3\ZZ)^{\oplus 2}$-extension. 
Therefore we have the rational maps $\mu_1$ and $\mu_2$ as desired.

\begin{rem}
It may be an interesting question to consider if there exists a  simple relation between $X_{1}$ and $X_{2}$ in the case of $A_{5}$ as above.
\end{rem}

Department of Mathematics and Information Science

Tokyo Metropolitan University

1-1 Minamiohsawa

 Hachoji 192-0397
 
 Tokyo
\end{document}